\documentclass[a4paper]{amsart}
\usepackage{lmodern}
\usepackage{amsthm}
\usepackage{amssymb}
\usepackage{amsmath}
\usepackage{hyperref}
\hypersetup{
  unicode=true,
  colorlinks=true,
  citecolor=blue,
  linkcolor=blue,
  anchorcolor=blue
}
\usepackage{todonotes}
\usepackage{tikz-cd}
\usepackage{dsfont,mathrsfs}
\usepackage{graphicx}

\usepackage[pagewise,mathlines]{lineno}

\usepackage[noadjust,nocompress]{cite}

\ifdefined\newaliascnt
  \relax
\else
  \usepackage{aliascnt}
\fi

\newtheorem{theorem}{Theorem}[section]
\newaliascnt{claim}{theorem}

\aliascntresetthe{claim}

\newaliascnt{lemma}{theorem}
\newtheorem{lemma}[lemma]{Lemma}
\aliascntresetthe{lemma}

\newaliascnt{proposition}{theorem}
\newtheorem{proposition}[proposition]{Proposition}
\aliascntresetthe{proposition}

\newaliascnt{corollary}{theorem}
\newtheorem{corollary}[corollary]{Corollary}
\aliascntresetthe{corollary}

\newaliascnt{fact}{theorem}
\newtheorem{fact}[fact]{Fact}
\aliascntresetthe{fact}

\theoremstyle{definition}
\newaliascnt{remark}{theorem}
\newtheorem{remark}[remark]{Remark}
\aliascntresetthe{remark}

\newaliascnt{definition}{theorem}
\newtheorem{definition}[definition]{Definition}
\aliascntresetthe{definition}

\newaliascnt{question}{theorem}
\newtheorem{question}[question]{Question}
\aliascntresetthe{question}

\newaliascnt{conjecture}{theorem}

\aliascntresetthe{conjecture}

\makeatletter
\newcommand{\getnamereftext}[1]{%
  \@ifundefined{r@#1}{}{%
    \unexpanded\expandafter
    \expandafter\expandafter\expandafter
    \expandafter\expandafter\expandafter
      {%
       \expandafter\expandafter\expandafter\@thirdoffive\csname r@#1\endcsname
      }%
  }%
}
\makeatother
\def\emptytemp{}
\newcommand{\futureref}[2]{
  \edef\therefname{\getnamereftext{#1}}
  \newtheorem*{futurethm}{\autoref{#1}}
  \ifx\therefname\emptytemp
    \begin{futurethm}
      #2
    \end{futurethm}
  \else
    \begin{futurethm}[\nameref*{#1}]
      #2
    \end{futurethm}
  \fi
  \let\futurethm\undefined
  \let\endfuturethm\undefined
}

\DeclareMathOperator{\rng}{rng}

\DeclareMathOperator{\id}{id}
\DeclareMathOperator{\crit}{crit}
\DeclareMathOperator{\Ult}{Ult}
\DeclareMathOperator{\Col}{Col}
\DeclareMathOperator{\otp}{otp}
\DeclareMathOperator{\cf}{cf}
\DeclareMathOperator{\acc}{acc}
\DeclareMathOperator{\tcl}{tcl}

\DeclareMathOperator{\sym}{sym}
\DeclareMathOperator{\fix}{fix}

\DeclareMathOperator{\Add}{Add}
\DeclareMathOperator{\supp}{supp}
\DeclareMathOperator{\meas}{\textsf{meas}}

\let\leq\leqslant
\let\nleq\nleqslant

\let\geq\geqslant

\newcommand{\ZF}{\mathsf{ZF}}
\newcommand{\ZFC}{\mathsf{ZFC}}

\newcommand{\SCH}{\mathsf{SCH}}
\newcommand{\DC}{\mathsf{DC}}
\newcommand{\HS}{\mathsf{HS}}

\newcommand{\AC}{\mathsf{AC}}

\newcommand{\WO}{\mathsf{WO}}
\newcommand{\SVC}{\mathsf{SVC}}
\newcommand{\KWP}{\mathsf{KWP}}
\newcommand{\ACH}{\mathsf{ACH}}
\newcommand{\Ord}{\mathrm{Ord}}
\newcommand{\HOD}{\mathrm{HOD}}
\newcommand{\HCD}{\mathrm{HCD}}

\newcommand{\LP}{\textup{\textsf{{\L}P}}}

\newcommand{\power}{\mathcal P}

\newcommand{\cW}{\mathcal W}

\newcommand{\cF}{\mathcal F}

\newcommand{\cS}{\mathcal S}
\newcommand{\cT}{\mathcal T}

\newcommand{\sG}{\mathscr G}
\newcommand{\sH}{\mathscr H}
\newcommand{\sF}{\mathscr F}
\newcommand{\sK}{\mathscr K}

\newcommand{\PP}{\mathbb P}

\newcommand{\RR}{\mathbb R}

\newcommand{\forces}{\mathrel{\Vdash}}

\newcommand{\1}{\mathds 1}

\newcommand{\tup}[1]{\langle#1\rangle}

\subjclass[2020]{Primary 03E55; Secondary 03E25, 03E35}
\keywords{supercompact cardinals, small violations of choice, HOD conjecture}
\date{November 3, 2026}
\title{Large cardinals, small violations of choice}
\author{Hope Duncan}
\address[H.~Duncan]{School of Mathematics, University of Leeds, Leeds, LS2~9JT, UK}
\email{mmhid@leeds.ac.uk}
\author{Gabriel Goldberg}
\address[G.~Goldberg]{Evans Hall, University Drive, Berkeley, CA 94720, United States}
\email{ggoldberg@berkeley.edu}
\author{Asaf Karagila}
\address[A.~Karagila]{School of Mathematics, University of Leeds, Leeds, LS2~9JT, UK}
\email{karagila@math.huji.ac.il}

\begin{document}
\begin{abstract}
We study the relationship between large cardinals and Blass' \emph{Small Violation of Choice} ($\SVC$), which states that the Axiom of Choice holds in a set-generic extension of the universe. Our main theorem is that in the presence of a supercompact cardinal, either the Axiom of Choice fails significantly, or else $\SVC$ holds in $L(\power(\Ord))$. We show that a supercompact is not strong enough to reflect $\SVC$, while an extendible cardinal is. If we assume Usuba's \emph{{\L}o\'s Property} holds for the supercompact, then we can outright derive $\SVC$ in the universe.
\end{abstract}
\maketitle
\section{Introduction}
In the past decade or so a great interest in the Axiom of Choice has been reignited in the set theoretic community. This is, in part, due to the works in \cite{Woodin:SEM1,Beyond} on large cardinal axioms which are incompatible with the Axiom of Choice. In general, large cardinal axioms extend the standard axioms of Zermelo--Fraenkel and the Axiom of Choice ($\ZFC$, or $\ZF$ when the Axiom of Choice is omitted), and they do so by positing the existence of some cardinal with additional properties from whence additional structure is added to the universe. These properties are often formulated as the existence of some non-trivial elementary embedding from the universe (or an initial segment thereof) into a transitive inner model, and the axioms form a natural hierarchy of strength by requiring this inner model is ``closer and closer'' to the full universe of sets in which we work.

Woodin showed that large cardinals provide us with additional structure in a very interesting way. Assuming a strongly compact cardinal exists,\footnote{Woodin's original assumptions were stronger, this one is due to Goldberg and was shown to be optimal.} the universe of sets is either ``very close'' to the inner model of all hereditarily ordinal definable sets ($\HOD$), in the sense that $\HOD$ agrees with the universe on the notion of singular cardinals and their successors, or else it is ``very far'' from it in that every regular cardinal in the universe is in fact a measurable cardinal in $\HOD$ (and more).

This dichotomy, known as \emph{The HOD Dichotomy}, suggested that even without the Axiom of Choice, large cardinal axioms produce a certain structure: we are either very close to an inner model in which the Axiom of Choice holds, or else we are hopelessly far from it.

We can understand ``close to an inner model where the Axiom of Choice holds'' in a different way using Blass' \cite{Blass:1979} principle of ``Small Violations of Choice'' ($\SVC$) which states that there is one set whose lack of well-order is the problem, when we wish to identify this set we will write $\SVC(S)$.

Blass showed that $\SVC(S)$ holds exactly when a generic well-ordering of $S$ will force the Axiom of Choice. Usuba \cite{Usuba:2021} showed that $\SVC$ implies that there is an inner model of the Axiom of Choice, over which the universe is a ``small extension''.

Throughout the paper we work in $\ZF$, unless stated otherwise. Our main theorem in this paper is the following dichotomy.
\futureref{thm:spc-dich}{Suppose that there exists a supercompact cardinal. Then exactly one of the following holds:
  \begin{enumerate}
  \item Every sufficiently large regular cardinal is $\omega$-J\'onsson.
  \item For every $\alpha\in\Ord$, $L(\power^\alpha(\Ord))\models\SVC$.
  \end{enumerate}
}

We show, therefore, that in the presence of large cardinals either the universe hopelessly fails the Axiom of Choice, or else the inner models $L(\power^\alpha(\Ord))$ are all very close to models where the Axiom of Choice holds. To that end, we prove a theorem that is essentially due to Woodin, stating that $\SVC$ propagates. Namely, if $L(\power(\Ord))\models\SVC$, then for every $\alpha$, $L(\power^\alpha(\Ord))\models\SVC$.

Our result builds on and improves Woodin's Theorem~229 in \cite{Woodin:SEM1} in a few ways: we replace its dependence on the HOD Conjecture with a precise dichotomy; we weaken its hypothesis to one supercompact cardinal; and we use the general technique of completely definable sets rather than Woodin's one-off construction.

On the question of whether (2) can be strengthened to $\SVC$ in $V$ itself is closely related to Woodin's \emph{Axiom of Choice Conjecture} \cite[Definition~3.51]{Woodin:HOD}, which states that $\SVC(V_\delta)$ holds when $\delta$ is an extendible cardinal. Both the HOD Conjecture and the Axiom of Choice Conjecture are sufficient to prove that there are no Reinhardt cardinals, and Woodin asked whether the HOD Conjecture implies the Axiom of Choice Conjecture.

We also answer a related question of Usuba, under the assumption of the \emph{{\L}o\'s Property} ($\LP$), which states that ultrapowers by measures on ordinals satisfy {\L}o\'s's Theorem.

\futureref{thm:spc+lp}{
Let $\kappa$ be a supercompact cardinal such that $\LP(\kappa)$ holds. Then $\SVC(V_{\kappa+1})$ holds.
}

Finally, we study the reflection of $\SVC$ below large cardinals. In particular, we show that if $\kappa$ is an extendible cardinal, then $\SVC$ implies $\SVC(V_\alpha)$ for some $\alpha<\kappa$. However, in the case of a supercompact cardinal we prove the following independence result:

\futureref{thm:spc+svc}{
  Let $\kappa$ be the least supercompact cardinal. There is a symmetric extension, $W$, in which $\kappa$ is supercompact, $\SVC(W_\kappa)$ holds, but $\SVC(W_\alpha)$ fails for all $\alpha<\kappa$.
}

To prove this theorem we also prove \autoref{thm:spc-lifting} which is a preservation theorem for supercompactness in symmetric extensions.

\subsection*{Acknowledgements}
The research was supported by UKRI Future Leaders Fellowship [MR/Y034058/1], by Engineering and Physical Sciences Research Council [UKRI5109], and by National Science Foundation under Grant No.\ 2401789. No data are associated with this article. For the purpose of open access, the authors have applied a CC-BY licence to any Author Accepted Manuscript version arising from this submission.

\subsection*{AI Declaration}
 The authors would like to apologise for any typos or misprints. AI was not used in the process of research or writing of this work.

\section{Preliminaries}
Let us begin with brief notational conventions. When $V$ is a model of $\ZF$ we write $V_\alpha$ to denote its $\alpha$th stage in the von Neumann hierarchy. When $X$ is a proper class, we write $L(X)$ to denote $\bigcup\{L(X\cap V_\alpha)\mid\alpha\in\Ord\}$, which we use most often with $L(\power^\alpha(\Ord))=\bigcup\{L(\power^\alpha(\kappa)\mid\kappa\in\Ord\}$. We write $\HOD$ to denote the class of all sets which are hereditarily ordinal definable. That is, $x\in\HOD$ if for every $a$ in the transitive closure of $\{x\}$, $a$ is definable using an ordinal parameter. We write $\HOD_X$ to allow parameters from the class $X$.

When $X$ is a set of ordinals and $\alpha$ is an ordinal, we write $[X]^\alpha$ to mean the set of subsets of $X$ whose order type is $\alpha$.\footnote{The notation is sometimes used to mean the set of subsets of \emph{cardinality} $\alpha$.}

For a set $X$, we write $\AC_X$ to denote the statement ``every family indexed by $X$ admits a choice function''. In this paper we will mostly focus on the case where $X$ is a well-ordered set, in which case we will simply consider it to be an well-ordered cardinal. We write $\AC_\WO$ to denote the statement ``every well-orderable family of sets admits a choice function''. For a study of various equivalents of $\AC_\WO$, see \cite[\S3]{RS:2024}.

Our treatment of forcing is standard (see, e.g.\ \cite{Jech:ST}).\footnote{We do not use class forcing in this work and so the term ``generic'' will always refer to \emph{set-generic}.} We want to remark that the Axiom of Choice is not generally used in the basic mechanism of forcing and genericity. It is, however, used in the proofs that chain conditions and closure properties of partial orders preserve cardinals. Nevertheless, variants of these without the Axiom of Choice can be shown quite easily.

\begin{definition}
  We say that $\SVC(X)$ holds if for every set $Y$ there exists an ordinal $\eta$ and a surjective map $f\colon X\times\eta\to Y$. In this case, we say that $X$ is a \emph{seed}. We write $\SVC$ to mean that there exists $X$ such that $\SVC(X)$ holds.
\end{definition}
The following fact is essentially proved in \cite[Theorem~4.6]{Blass:1979}.
\begin{fact}\label{fact:SVC}
  For any set $X$, $\SVC(X^{<\omega})$ holds if and only if $\Col(\omega,X)=X^{<\omega}$ forces the Axiom of Choice.
\end{fact}
For two sets $A$ and $B$ we write $A\leq B$ if there is an injection from $A$ into $B$ and $A\leq^*B$ if there is a partial surjection from $B$ onto $A$.\footnote{We need the surjection to be partial to account for the empty set.} For a set $X$ we let $\aleph(X)$ be $\min\{\alpha\in\Ord\mid\alpha\nleq X\}$ and $\aleph^*(X)=\min\{\alpha\in\Ord\mid\alpha\nleq^*X\}$.\footnote{These are known as the \emph{Hartogs} and \emph{Lindenbaum} numbers of $X$.}
\begin{lemma}\label{lemma:svc-ordinal}
  Suppose that $f\colon X\times\eta\to Y$ is a surjective map and suppose that $\eta$ is the least ordinal for which such a map exists. Then $\eta<\aleph^*(Y)$.
\end{lemma}
\begin{proof}
  We let $g\colon Y\to\eta$ be $g(y)=\min\{\alpha<\eta\mid\exists x\in X,f(x,\alpha)=y\}$. Easily, $g$ is well-defined. Moreover, $f\restriction X\times\rng(g)$ is still a surjection onto $Y$, and therefore by the minimality of $\eta$, it must be that $\otp(\rng g)=\eta$.
\end{proof}
Let us denote by $\cW(a)$ the set of all well-ordered sequences of elements of $a$. That is, the set of all injective functions from an ordinal into $a$. This is a set, as the domain of any member of $\cW(a)$ must be below $\aleph(a)$.
\begin{lemma}\label{lemma:htg-increase}
  For every set $X$ and ordinal $\eta$, $X^\eta\leq\cW(X)\times\aleph(X)^\eta$.
\end{lemma}
\begin{proof}
  If $f\in X^\eta$, we can recursively define some $A\subseteq\eta$ on which $f$ is injective and has the same range. Let $\alpha=\otp(A)$ and $\pi\colon\eta\to\alpha$ an extension of the Mostowski collapse of $A$. Then $f\mapsto(\pi\circ f\restriction A,\pi)$ is an injective function as wanted.
\end{proof}
We shall prove the following propagation theorem for $\SVC$. This is similar to Woodin's claim \cite[Remark~232]{Woodin:SEM1}.
\begin{theorem}\label{thm:svc-prop}
  If $L(\power(\Ord))\models\SVC$, then for every $\alpha$, $L(\power^\alpha(\Ord))\models\SVC$.
\end{theorem}
To prove the theorem we first prove a lemma for the successor steps.
\begin{lemma}\label{lemma:svc-prop-succ}
  Let $X$ be a transitive class and $\power(\Ord)\subseteq X$. If $L(X)\models\SVC(S)$, then $L(\power(X))\models\SVC(\cW(\power(S))\times\power(S))$.
\end{lemma}
\begin{proof}
  Let $S$ be a seed for $\SVC$ in $L(X)$ and work inside the model $L(\power(X))$. We let $\gamma=\aleph(\power(S))$.

  For some $x\in X$, let $\eta$ be some ordinal such that $x\leq^*\eta\times S$ in $L(X)$. Then $\power(x)\leq^*\power(\eta\times S)$,\footnote{This is far from optimal, of course.} as $\power(\eta\times S)$ has the same cardinality as $\power(S)^\eta$, \autoref{lemma:htg-increase} allows us to conclude that $\power(S)^\eta\leq^*\cW(\power(S))\times\gamma^\eta$.

  As $L(X)$ computes $\gamma^\eta$ correctly, there is some $\delta$ such that $\gamma^\eta\leq^* S\times\delta$. Thus, $\power(x)\leq^*\cW(\power(S))\times\power(S)\times\delta$. Now $\SVC(\cW(\power(S))\times\power(S))$ holds in $L(\power(X))$.
\end{proof}
\begin{proof}[Proof of \autoref{thm:svc-prop}]
  We construct a sequence of seeds, $S_\alpha$, for $L(\power^\alpha(\Ord))$ for $\alpha>0$ from a fixed seed, $S_1$, for $L(\power(\Ord))$, which can be taken as $\power(\delta)$ for some cardinal $\delta$.

  Suppose $S_\alpha$ is a seed for $L(\power^\alpha(\Ord))$. Letting $S_{\alpha+1}=\cW(\power(S_\alpha))\times\power(S_\alpha)$, as computed in $L(\power^{\alpha+1}(\Ord))$, then by \autoref{lemma:svc-prop-succ}, $S_{\alpha+1}$ is a seed.

  Suppose $\alpha$ is a limit ordinal and for all $\beta<\alpha$, $S_\beta$ was defined. We define \[Y=\bigcup\{\{\beta\}\times S_\beta^{<\omega}\mid\beta<\alpha\},\] and set $S_\alpha=\cW(\power(Y))$. We claim that $S_\alpha$ is a seed in $L(\power^\alpha(\Ord))$.

  It is enough to show that for every $\kappa$ there is some $\eta$ such that $\power^\alpha(\kappa)\leq^*S_\alpha\times\eta$. First, note that for each $\beta<\alpha$, there is some $\eta_\beta$ such that $\power^\beta(\kappa)\leq^*S_\beta\times\eta_\beta$. Let $\eta=\sup\{\eta_\beta\mid\beta<\alpha\}$, then every $\power^\beta(\kappa)$ is the transitive collapse of some $D\in\power(S_\beta^{<\omega}\times\eta)$. In particular, there is a partial surjection from $\power(S_\alpha)^\eta$ onto $\power^\alpha(\kappa)$, by taking the Mostowski collapse of any set encoding an element of $\power^\alpha(\kappa)$.

  By the fact that $\aleph(\power(Y))^\eta$ is computed correctly in $L(\power(\Ord))$, we have that \[\power(Y)^\eta\leq^*S_\alpha\times\aleph(\power(Y))^\eta\leq^*S_\alpha\times S_1\times\delta.\] It is not hard to see that $S_1\times S_\alpha$ is bijective with $S_\alpha$, and so we have $\power^\alpha(\delta)\leq^*S_\alpha\times\delta$ as wanted.
\end{proof}
\begin{definition}
  We say that $\alpha$ is a \emph{$\beth$-fixed point} if for every $x\in V_\alpha$, $\aleph^*(x)<\alpha$.
\end{definition}
The following is an easy observation about $\beth$-fixed points.
\begin{fact}
  For an ordinal $\alpha$ the following are equivalent.
  \begin{enumerate}
  \item Every well-ordered relation in $V_\alpha$ is isomorphic to an ordinal in $V_\alpha$.
  \item $V_\alpha$ satisfies the Mostowski collapse lemma.
  \item $\alpha$ is a $\beth$-fixed point.
  \end{enumerate}
\end{fact}
\begin{lemma}\label{lemma:strong-limit-svc}
  Suppose that $\SVC(X)$ holds. If $\alpha$ is a $\beth$-fixed point and $X\in V_\alpha$, then $V_\alpha\models\SVC(X)$.
\end{lemma}
\begin{proof}
  If $a\in V_\alpha$, then there is some $\eta$ and a surjection $f\colon\eta\times X\to a$. By \autoref{lemma:svc-ordinal}, the least such $\eta$ is below $\aleph^*(a)<\alpha$.
\end{proof}
\begin{definition}
  For $\alpha\in\Ord$ we say that the \emph{Kinna--Wagner* $\alpha$ Principle} holds, if for every set $X$ there an ordinal $\eta$ and a surjection from $\power^\alpha(\eta)$ onto $X$. We denote this stated by $\KWP_\alpha^*$.
\end{definition}
This principle was recently explored by Karagila and Schilhan in \cite{KS:KWP} where the following facts can be found.
\begin{fact}\label{fact:KWP}The following facts are proved in \cite{KS:KWP}.
  \begin{enumerate}
  \item $L(\power^\alpha(\Ord))\models\KWP_\alpha^*$. (Proposition~3.9)
  \item If $V\models\ZFC$ and $V[G]$ is a generic extension of $M$, then for every intermediate model of $\ZF$, $V\subseteq M\subseteq V[G]$, $M\models\KWP_\alpha^*$ for some $\alpha$ if and only if $M\models\SVC$.\footnote{This is the particular case of the theorem when $V$ satisfies $\AC$.} (Theorem~5.5)
  \end{enumerate}
\end{fact}
\section{A dichotomy theorem for the Axiom of Choice Conjecture}
Large cardinals have a significant effect on the structure of the universe. Woodin proved the \emph{HOD Dichotomy} \cite{Woodin:HOD} which is analogous to Jensen's covering lemma for the inner model $\HOD$. It states that in the presence of an extendible cardinal, the universe of sets is very close to $\HOD$ or very far from $\HOD$. This was later improved by Goldberg to a strongly compact cardinal \cite{Goldberg:strcompt}, and Ben-Neria and Hayut \cite{BNH} showed that the large cardinal assumptions are necessary. The \emph{HOD Hypothesis} is the assertion stating that the universe is close to $\HOD$. Here we study another conjecture of Woodin, named the \emph{Axiom of Choice Conjecture}, which states ``if $\delta$ is an extendible cardinal, then $\SVC(V_\delta)$ holds''.

We will then show that if a supercompact cardinal exists, then there is a slightly weaker dichotomy at play: either almost all regular cardinals are $\omega$-J\'onsson, or $L(\power^\alpha(\Ord))$ satisfies $\SVC$.\footnote{If we assume a proper class of supercompact cardinals, this can be phrased as either the \emph{HOD Conjecture}\footnotemark\ fails or $\SVC$ holds in these inner models.}\footnotetext{The HOD Conjecture states that $\ZFC$ proves the HOD Hypothesis.}

We then show that under Usuba's \emph{{\L}o\'s's Property}, $\SVC$ holds in the universe itself. In the next section we show that supercompactness is not sufficient to reflect $\SVC$ like extendibility.

\subsection{A dichotomy theorem for a supercompact cardinal}
\begin{definition}
  We say that $\kappa$ is a \emph{$\lambda$-supercompact cardinal} if there exists $\beta\geq\lambda$ and a transitive set $N$ such that $N^{V_\lambda}\subseteq N$, and there is an elementary embedding $j\colon V_\beta\to N$ with $\crit(j)=\kappa$ and $j(\kappa)>\lambda$. If $\kappa$ is $\lambda$-supercompact for all $\lambda$ we say that $\kappa$ is a supercompact cardinal.\footnote{Usuba refers to this as \emph{Woodin-supercompactness}, reserving the unqualified supercompactness to the second-order definition where $j$ is defined on $V$.}
\end{definition}

We will use the term \emph{$\lambda$-supercompact embedding} to mean that an embedding witnessing $\lambda$-supercompactness. Namely, $j\colon V_\alpha\to N$ with $\alpha\geq\lambda$ and $N^{V_\lambda}\subseteq N$. In the standard context of $\ZFC$ we would use this term to mean closure under $\lambda$-sequences, however in all of our applications where the Axiom of Choice holds we will have $\lambda=|V_\lambda|$, so no confusion will arise.

\begin{definition}
  We say that a cardinal $\lambda$ is \emph{$\omega$-J\'onsson} if for every $f\colon[\lambda]^\omega\to\lambda$, there exists $H\in[\lambda]^\lambda$ such that $f``[H]^\omega\neq\lambda$.
\end{definition}
By \cite{ErdosHajnal}, the Axiom of Choice implies that no $\omega$-J\'onsson cardinals exist.

\begin{theorem}[Axiom of Choice Dichotomy]\label{thm:spc-dich}
  Suppose that there exists a supercompact cardinal. Then exactly one of the following holds:
  \begin{enumerate}
  \item Every sufficiently large regular cardinal is $\omega$-J\'onsson.
  \item For every $\alpha\in\Ord$, $L(\power^\alpha(\Ord))\models\SVC$.
  \end{enumerate}
\end{theorem}
Following Woodin's terminology, we shall refer to the failure of (1) in the above dichotomy as ``\emph{the Axiom of Choice Hypothesis}'',\footnote{Technically, this should be called the \emph{weak} Axiom of Choice Hypothesis, as $\SVC$ only holds in an inner model of $V$.} abbreviated as $\ACH$ or $\ACH(\kappa)$ when specifying a supercompact cardinal, $\kappa$. Our goal, therefore, is to show that the Axiom of Choice Hypothesis implies $\SVC$ holds in $L(\power^\alpha(\Ord))$. And by \autoref{thm:svc-prop}, it is enough to argue $\SVC$ holds in $L(\power(\Ord))$.

Let us first show the Axiom of Choice Hypothesis holds when $L(\power(\Ord))\models\SVC$.

\begin{proposition}\label{prop:oJ-solovay}
  If $\lambda$ is a regular cardinal and it is $\omega$-J\'onsson, then $S^\lambda_\omega$ cannot be partitioned into $\lambda$-many disjoint stationary sets.
\end{proposition}
\begin{proof}
  Suppose that $\{S_\alpha\mid\alpha<\lambda\}$ is a partition of $S^\lambda_\omega$ into stationary sets. We define $f\colon[\lambda]^\omega\to\lambda$ by $f(x)=\min\{\alpha\mid\sup x\in S_\alpha\}$. For every $S\in[\lambda]^\lambda$ note that the set of accumulation points of $S$, that is, $\acc(S)=\{\beta<\lambda\mid\beta=\sup(\beta\cap S)\}$, is a club in $\lambda$. Therefore, for any $\alpha<\lambda$ we may find some $\beta\in S_\alpha\cap\acc(S)$, which means that there is some $s\in [S\cap\beta]^\omega$ such that $\sup s=\beta\in S_\alpha$.

  Therefore, $\lambda$ cannot be $\omega$-J\'onsson if $S^\lambda_\omega$ can be partitioned into $\lambda$-many disjoint stationary sets.
\end{proof}

\begin{theorem}
  If $L(\power(\Ord))\models\SVC$, then for a proper class of cardinals, $\lambda$, $\lambda^+$ is regular and not $\omega$-J\'onsson.
\end{theorem}
\begin{proof}
  Let $N=L(\power(\Ord))$, and suppose that $N\models\SVC$. Then there is some $\PP\in N$ and $G\subseteq\PP$ which is $N$-generic filter, such that $N[G]\models\ZFC$. If follows that for all $\lambda>|\PP|^{N[G]}$, $\lambda^+$ is regular, so $\lambda^+$ must be regular in $N$.

  Moreover, $(S^{\lambda^+}_\omega)^N$ remains stationary: otherwise $C$ is a club disjoint from it and some $p\in G$ have decided a club subset of $C$ which is already in $N$. Similarly, if we partition $(S^{\lambda^+}_\omega)^N$ in $N[G]$ we may pull the partition back to $N$: some $p\in G$ will have decided $\lambda^+$-many stationary sets.

  However, $N$ and $V$ agree on sets of ordinals, including notions such as regularity and stationary sets. Therefore, by \autoref{prop:oJ-solovay}, $\lambda^+$ is not $\omega$-J\'onsson.
\end{proof}

We shall now proceed to prove that $\ACH$ implies $\SVC$ in $L(\power(\Ord))$.

\begin{definition}
  For a set $X$ and a well-ordered cardinal $\kappa$, $\meas_\kappa(X)$ denotes the set of all $\kappa$-complete measures on $X$. We write $\meas_\kappa$ to denote $\bigcup_{\delta\in\Ord}\meas_\kappa(\delta)$.
\end{definition}

We will need the following fact, due to Goldberg and stated here in the special case relevant to our proof, in order to prove that $\meas_\kappa(\delta)$ is well-orderable when $\kappa$ is a supercompact.

\begin{fact}[{\cite[Lemma~3.5]{Goldberg:2024meas}}]\label{fact:completeness}
  Let $\kappa$ be a $\beth$-fixed point. For every $U\in\meas_\kappa$ and $I\in V_\kappa$, if $\{A_i\mid i\in I\}\subseteq U$, then $\bigcap_{i\in I}A_i\in U$.
\end{fact}

\begin{lemma}[Kunen]\label{lemma:meas-from-ordinals}
  Let $\kappa$ be a supercompact cardinal and let $\delta\geq\kappa$. For every $\eta>\delta+2$ and an $\eta$-supercompact embedding, $j\colon V_\eta\to N$, if $U\in\meas_\kappa(\delta)$, there is some ordinal $\alpha<j(\delta)$ such that $A\in U$ if and only if $\alpha\in j(A)$.
\end{lemma}
\begin{proof}
  Let $j\colon V_\eta\to N$ be an $\eta$-supercompact embedding with $\eta>\delta+2$. Since $U\in V_\eta$ and $N^{V_\eta}\subseteq N$, $j``U\in N$. As $j(\kappa)$ is a regular $\beth$-fixed point in $N$ and $j(U)$ is $j(\kappa)$-complete, we have that $\bigcap j``U\in j(U)$. In particular, $\bigcap j``U\neq\varnothing$. Let $\alpha=\min\bigcap j``U$, then $A\in U$ if and only if $\alpha\in j(A)$ as wanted.
\end{proof}

\begin{corollary}\label{lemma:spc-wo-meas}
  If $\kappa$ is supercompact, then $\meas_\kappa(\delta)$ is a well-orderable set.
\end{corollary}
\begin{proof}
  Fix an embedding as in \autoref{lemma:meas-from-ordinals} and define \[U_0<_j U_1\iff\min\left(\bigcap j``U_0\right)<\min\left(\bigcap j``U_1\right).\qedhere\]
\end{proof}

\begin{lemma}[Solovay]\label{lemma:spc-wo-con}
  If $\kappa$ is a supercompact cardinal and $\delta>\kappa$ is a regular cardinal that is not $\omega$-J\'onsson, then there is a supercompact measure on $\power_\kappa(\delta)$ concentrating on a well-orderable set.
\end{lemma}
\begin{proof}
  Let $j\colon V_\alpha\to N$ be a $\delta$-supercompact embedding with $\alpha>\delta+2$. Since $\delta$ is not $\omega$-J\'onsson, fix a witnessing function $f\colon[\delta]^\omega\to\delta$ for which whenever $A\in[\delta]^\delta$, $f``[A]^\omega=\delta$. Finally, let $U$ be the supercompact measure derived from $j$ and $j``\delta$,\[U=\{X\subseteq\power_\kappa(\delta)\mid j``\delta\in j(X)\}.\]

  For $\eta<\delta$, let $\sigma_\eta$ be the least $\omega$-club, if such a club exists, which is closed under $f$. Namely, $f``[\sigma_\eta]^\omega\subseteq\sigma_\eta$, and it is the least $\omega$-club of $\eta$ with this property.

  We claim that the set $A=\{\sigma_\eta\mid \eta<\delta,\sigma_\eta\text{ exists, and }|\sigma_\eta|<\kappa\}$ belongs to $U$. As the supremum of $\sigma_\eta$ is unique, the set $A$ can be well-ordered, which will complete our proof.

  Let us check that $j``\delta\in j(A)$. That is to say, letting $\eta=\sup j``\delta$, we need to check that $\sigma_\eta=j``\delta$.

  Since $j$ is a supercompact embedding, $j``\delta$ is an $\omega$-club in $\eta$ and $j(f)``[j``\delta]^\omega=j``\delta$. It remains to show that it is in fact the least such $\omega$-club.

  Suppose that $S\subseteq\eta$ was an $\omega$-club such that $j(f)``[S]^\omega\subseteq S$. We let $T$ be the preimage of $S$, that is, $T=\{\alpha<\delta\mid j(\alpha)\in S\}$. Then $T$ must be unbounded in $\delta$, as $S\cap j``\delta$ is an $\omega$-club in $\eta$. By the regularity of $\delta$, $T\in[\delta]^\delta$, and so $f``[T]^\omega=\delta$. In particular, $j(f)``[S]^\omega\subseteq j(f)``[j``T]^\omega=j``\delta$, so $j``\delta\subseteq S$ as wanted.
\end{proof}

\begin{definition}
  We say that $X$ is a \emph{$\kappa$-completely definable set} if it is definable from a parameter in $\meas_\kappa$. We denote by $\mathsf{CD}_\kappa$ the class of $\kappa$-completely definable sets. We say that $X$ is \emph{hereditarily $\kappa$-completely definable} if $\tcl(\{X\})\subseteq\mathsf{CD}_\kappa$ and we denote the class of hereditarily $\kappa$-definable sets by $\HCD_\kappa$.

  We will use the notation $\HCD_{\kappa,S}$ to denote the relativised definition. Namely, the \emph{elements of $S$} can be used as parameters in the definition.
\end{definition}

As the ordinals correspond to the principal ultrafilters, and finite sequences of ultrafilters correspond to their Fubini product, it is easy to see that $\HOD\subseteq\HCD_\kappa$. Indeed, $\HCD_\kappa=\HOD_{\meas_\kappa}$.

\begin{lemma}\label{lemma:ord-hcd-capture}
  Let $\kappa$ be a supercompact cardinal and let $\delta>\kappa$ be regular which is not $\omega$-J\'onsson. Then for every $A\subseteq\delta$ there is $a\subseteq\kappa$ such that $A\in\HCD_{\kappa,\{a\}}$.
\end{lemma}
\begin{proof}
  Let $U'$ be a supercompact measure on $\power_\kappa(\delta)$ given by \autoref{lemma:spc-wo-con}, let $\eta$ be the order type of the well-orderable set constructed in said proof, and let $U$ be the $\kappa$-complete ultrafilter on $\eta$ given by the obvious projection of $U'$ onto $\eta$.

  Consider the structures $\Ult(\power_\kappa(\kappa),U)\cong\Ult(\power_\kappa(\kappa),U')$, let us first notice that the map $x\mapsto\otp(x)$ is a function $\power_\kappa(\lambda)\to\kappa$ represents $\delta$ in $\Ult(\power_\kappa(\kappa),U')$. So, there is some $g\colon\eta\to\power_\kappa(\kappa)$ which represents $\delta$.

  Therefore, given any $A\subseteq\delta$, the function $x\mapsto\otp(A\cap x)$ represents $A$ in the ultrapower, and so there is some $g_A\colon\eta\to\power_\kappa(\kappa)$ such that $[g_A]_U=A$.

  In particular $A$ is definable from $U$ and $\tup{g_A(\xi)\mid\xi<\eta}$. Since $g_A(\xi)\in\power_\kappa(\kappa)$, $g_A$ can be seen as a subset of $\kappa\times\kappa\times\kappa$,\footnote{E.g., $(\xi,\eta,\alpha)\in g_A$ suggests that $\alpha\in g_A(\xi)\subseteq\eta$.} and so coded by $a\subseteq\kappa$. Therefore, $A\in\HCD_{\kappa,\{a\}}$.
\end{proof}
\begin{corollary}
  If $\ACH(\kappa)$ holds, then $L(\power(\Ord))\subseteq\HCD_{\kappa,\power(\kappa)}$.
\end{corollary}
\begin{proof}
  For every regular $\delta$ that is not $\omega$-J\'onsson, $\power(\delta)\subseteq\HCD_{\kappa,\power(\kappa)}$. As there is a proper class of such $\delta$, every set of ordinals is in $\HCD_{\kappa,\power(\kappa)}$.
\end{proof}
\begin{lemma}\label{lemma:hcd-zfc}
  Assuming $\ACH(\kappa)$ there exists $a\subseteq\kappa$ such that $\HCD_{\kappa,\{a\}}\models\ZFC$.
\end{lemma}
\begin{proof}
  For any $\alpha$, there exists some large enough $\delta$ such that every element of $V_\alpha\cap\HCD_\kappa$ can be defined from a measure in $\meas_\kappa(\delta)$. Fixing a well-ordering of $\meas_\kappa(\delta)$ induces a well-ordering of $V_\alpha\cap\HCD_\kappa$. Therefore, there is a least $\eta_\alpha$ for which a subset of $\eta_\alpha$ codes the structure $V_\alpha\cap\HCD_\kappa$. By \autoref{lemma:ord-hcd-capture}, there is some $a\subseteq\kappa$ such that such a set of ordinals lies in $\HCD_{\kappa,\{a\}}$.

  Setting $A_\alpha$ as the set of all such $a$, there is a proper class of ordinals for which $A_\alpha=A$. Pick any $a\in A$, and we have that every set in $\HCD_\kappa$ can be well-ordered inside $\HCD_{\kappa,\{a\}}$.

  Finally, if we can show that $\meas_\kappa\cap V_\alpha$ has a $\kappa$-complete definable well-ordering, we can apply Vop\v{e}nka's theorem to get that $a$ lies in a generic extension of $\HCD_\kappa$. It follows, then, that $\HCD_{\kappa,\{a\}}$ is a generic extension of $\HCD_\kappa$. In particular, every set in $\HCD_{\kappa,\{a\}}$ is the image of a set in $\HCD_\kappa$ which can be well-ordered there, and so $\HCD_{\kappa,\{a\}}\models\ZFC$.

  To see that $\meas_\kappa\cap V_\alpha$ has a $\kappa$-complete definable well-ordering first note that \autoref{lemma:spc-wo-meas} implies that for a large enough $\delta$ we can define such well-ordering from a supercompact measure on $\power_\kappa(\delta)$, and by \autoref{lemma:spc-wo-con}, such measure is isomorphic to one defined on an ordinal, as wanted.
\end{proof}
\begin{corollary}\label{cor:ach-hcd-svc}
  If $\ACH(\kappa)$ holds, then $\HCD_{\kappa,\power(\kappa)}\models\SVC(\power(\kappa))$.
\end{corollary}
\begin{proof}
  Fix $a\subseteq\kappa$ such that $\HCD_{\kappa,\{a\}}\models\ZFC$. Given $X\in\HCD_{\kappa,\power(\kappa)}$, each $x\in X$ is definable over $\HCD_{\kappa,\{a\}}$ using some $b\subseteq\kappa$. Let $\delta$ be large enough such that $\meas_\kappa(\delta)$ contains all the necessary ultrafilters for defining all the members of $X$ in some $V_\alpha$. As $\meas_\kappa(\delta)\in\HCD_\kappa$, it can be well-ordered in $\HCD_{\kappa,\{a\}}$. We fix such a well-ordering, $\{U_\alpha\mid\alpha<\eta\}$, an enumeration of the formulas in the language of set theory, $\{\varphi_n\mid n<\omega\}$, and a coding function $\pi\colon\power(\kappa)\to\power(\kappa)^{<\omega}$. We define the partial surjection \[F(\omega\cdot\alpha+n,b) = x\iff \varphi_n(x,U_\alpha,\pi^{-1}(b))\text{ defines }x\text{ in }V_\alpha.\]

  Then $F$ is internal to $\HCD_{\kappa,\power(\kappa)}$ and shows that $\SVC(\power(\kappa))$ holds as wanted.
\end{proof}
\begin{proof}[Proof of \autoref{thm:spc-dich}]
  It remains to prove that under $\ACH(\kappa)$, $L(\power(\Ord))\models\SVC$. Towards that end, we will find an inner model of $\ZFC$, $M$, such that for some $M$-generic filter, $G$, $M\subseteq L(\power(\Ord))\subseteq M[G]$.

  Then, by \autoref{fact:KWP}(1), $L(\power(\Ord))\models\KWP_1^*$ and $M\models\ZFC$. Therefore, by \autoref{fact:KWP}(2), $L(\power(\Ord))^V\models\SVC$ as wanted.

  We let $M=\HCD_{\kappa,\{a\}}$ in which $\ZFC$ holds, guaranteed to exist by \autoref{lemma:hcd-zfc} and $N=\HCD_{\kappa,\power(\kappa)}$. Note that $M=\HOD_M=\HOD_{M\cup\{\power(\kappa)\}}$ and $N=\HOD_N$, as $\power(\kappa)$ is itself ordinal definable. By Theorem~3.30 of \cite{GS:Part1}, there exists a $V$-generic filter $g\subseteq\Col(\omega,\power(\kappa))$ which is contained in an $M$-generic extension and $N$-generic extension. In particular, $N[g]$ satisfies the Axiom of Choice. We verify Bukovsk\'y's condition \cite{Bukovsky:1973,Bukovsky:2017} for $N[g]$ to be a generic extension of $M$, which will then provide us with the models necessary for the theorem.

  Let $f\in N[g]$ be a function $f\colon\kappa\to\Ord$, and let $\dot f\in N$ be a $\PP_N$-name for this function. Then there exists some $b\subseteq\kappa$ such that $\dot f$ is definable by $\varphi(x,U,b)$ for some measure $U$. Varying over all $c\subseteq\kappa$, we let $\dot f_c$ denote the name defined by $c\subseteq\kappa$, if it is a name at all. Letting \[F(\alpha) = \{\beta\mid(\exists p\in\Col(\omega,\power(\kappa)))(\exists c\subseteq\kappa)(p\forces\dot f_c(\check\alpha)=\check\beta\},\] we immediately have that $F\in M$ and it is not hard to check that $f(\alpha)\in F(\alpha)$ for all $\alpha$, and $|F(\alpha)|<\aleph^*(\power(\kappa))$.
\end{proof}
\subsection{Supercompactness and the {\L}o\'s Property}
\begin{definition}
  Let $I$ be a set and let $U$ be an ultrafilter over $I$. We say that $U$-$\AC_I$ holds if whenever $\{A_i\mid i\in I\}$ is a family of non-empty sets, then there is some $X\in U$ such that $\prod_{i\in X}A_i\neq\varnothing$.
\end{definition}
Mitchell Spector showed in \cite{Spector:1988} that $U$-$\AC_I$ is equivalent to the statement that the ultrapower embedding $j\colon V\to V^I/U$ is elementary, as well as that {\L}o\'s's theorem holds for ultraproducts by $U$, etc. This motivated Toshimichi Usuba to formulate the following property.
\begin{definition}
  We say that the \emph{{\L}o\'s Property holds at $\kappa$}, abbreviated as $\LP(\kappa)$, if for every $\lambda\geq\kappa$ and every $\kappa$-complete ultrafilter on $\lambda$, $U$-$\AC_\lambda$ holds.\footnote{Usuba explicitly requires that the ultrafilter is not $\kappa^+$-complete, but all of our applications are in the case where $\kappa$ is supercompact and the ultrafilter is derived from an embedding whose critical point is $\kappa$.}
\end{definition}

As an easy observation, $\AC_\kappa$ implies $U$-$\AC_\kappa$ for any ultrafilter $U$ over $\kappa$. Therefore, $\AC_\WO$ implies that $\LP(\kappa)$ holds for all $\kappa$. It is consistent, however, that $\LP(\kappa)$ holds for all $\kappa$, while $\AC_\WO$ fails. For example, if one construct Cohen's first model over any model of $\ZFC$ containing large cardinals, then $U$-$\AC_\kappa$ holds for any uncountable $\kappa$, but $\AC_\WO$ fails (in fact, $\AC_\omega$ fails).

The following bit of notation will be helpful in this section. Given an elementary embedding $j\colon V_\alpha\to M$, we let $D(j,\delta)$ be the derived measure from $\delta$, given by $\{A\subseteq\delta\mid\sup j``\delta\in j(A)\}$. It is not hard to check that $D(j,\delta)$ is $\crit(j)$-complete, and if $M^\delta\subseteq M$, then it is a uniform ultrafilter which extends the club filter on $\delta$.

\begin{theorem}\label{thm:spc+lp}
  Let $\kappa$ be a supercompact cardinal such that $\LP(\kappa)$ holds. Then $\SVC(V_{\kappa+1})$ holds.
\end{theorem}

For the remainder of the section, let $\kappa$ be a fixed supercompact cardinal for which $\LP(\kappa)$ holds. We begin by exploring the structure consequences of this axiom before proving the theorem.

\begin{lemma}\label{lemma:lp-reg}
  For all $\delta\geq\kappa$, $\delta^+$ is regular.
\end{lemma}
\begin{proof}
  Let $j\colon V_\eta\to M$ be a $\delta^+$-supercompact embedding with critical point $\kappa$. Let $D=D(j,\delta^+)$, then $D$ is a $\kappa$-complete ultrafilter on $\delta^+$. Moreover, if $A\in D$, then $\sup A=\delta^+$, otherwise let $\alpha=\sup A<\delta^+$, then $j(\alpha)<\sup j``\delta^+<j(\alpha)$ which is impossible.

  For each $\alpha<\delta^+$ let $A_\alpha=\{f\colon\alpha\to\delta\mid f\text{ is injective}\}$. By $\LP(\kappa)$, there exists a set $A\in D$ and a function defined on $A$ such that $f(\alpha)\in A_\alpha$ for all $\alpha\in A$. Therefore, we can uniformly enumerate all the $\alpha<\delta^+$, and so it must be that $\cf(\delta^+)=\delta^+$. Otherwise $\delta^+$ is a union of uniformly enumerated sets of size $\delta$.
\end{proof}
\begin{lemma}\label{lemma:fodor}
  For every regular $\delta\geq\kappa$, the club filter on $\delta$ is normal.
\end{lemma}
\begin{proof}
  Let $\delta\geq\kappa$ be a regular cardinal and let $\cF$ denote the club filter on $\delta$. Since $\kappa$ is a supercompact cardinal, $\cF$ is $\kappa$-complete. Let $j\colon V_\alpha\to M$ be a $\delta$-supercompact embedding.

  Let $\gamma\leq\delta$, and suppose that $\cF$ is $\gamma$-complete. Let $\{A_\alpha\mid\alpha<\gamma\}\subseteq\cF$, which we may assume by the completeness, is a decreasing sequence of sets. It suffices to show that there is a sequence of clubs, $C_\alpha\subseteq A_\alpha$, as the intersection (or diagonal intersection in the case $\gamma=\delta$) of clubs is always a club. Moreover, as the sequence is decreasing, it suffices to find such sequence of clubs on any cofinal subset of $\gamma$.

  Towards that goal, let $D=D(j,\gamma)$, then $D$ is a $\kappa$-complete and uniform ultrafilter on $\gamma$, and so by $\LP(\kappa)$, $D$-$\AC_\gamma$ holds, so we can choose clubs as wanted, and the conclusion follows.
\end{proof}

\begin{corollary}
  For every regular cardinal $\delta\geq\kappa$ and stationary $S\subseteq\delta$, there exists a stationary partition of $S$ into $\delta$ parts, $\tup{S_\alpha\mid\alpha<\delta}$.
\end{corollary}
\begin{proof}
  As the club filter on $\delta$ is normal, Fodor's lemma holds, and so the usual proof of Solovay's theorem can be carried out. When we need to choose clubs, we may fix some $\kappa$-complete ultrafilter $D$ on $\delta$ and use $\LP(\kappa)$ to select them on a cofinal set in $\delta$.
\end{proof}
Using \autoref{prop:oJ-solovay} we immediately obtain the following corollary.
\begin{corollary}
  $\ACH(\kappa)$.\qed
\end{corollary}

\begin{lemma}\label{lemma:stationary-coding}
  For every regular cardinal $\delta\geq\kappa$ and every stationary partition of $S^\delta_\omega$, $\vec S=\tup{S_\alpha\mid\alpha<\delta}$, if $j$ is a $\delta$-supercompact embedding, then $j(\vec S)_\alpha$ reflects at $\rho=\sup j``\delta$ if and only if $\alpha=j(\beta)$ for $\beta<\delta$.
\end{lemma}
\begin{proof}
  Since $j``S^\delta_\omega$ is an $\omega$-club in $\rho$, $j(\vec S)_\alpha$ reflects at $\rho$ if and only if $j(\vec S)_\alpha\cap j``S^\delta_\omega$ is a stationary set. However, $j``S_\beta\subseteq j(S_\beta)\cap\rho$, and $\{j``S_\alpha\mid\alpha<\delta\}$ is a partition of $j``S^\delta_\omega$. So, if $\alpha\notin j``\delta$, $j(\vec S)_\alpha\cap j(S_\beta)=\varnothing$. Therefore, $j(\vec S)_\alpha\cap\rho\cap j``S^\delta_\omega=\varnothing$ as well. The other direction is trivial as $j$ is $\delta$-supercompact and $j``\delta$ is in the target model.
\end{proof}

\begin{lemma}\label{lemma:it-capturing}
  For every $\alpha$ there is some $U\in\meas_\kappa$ such that $V_\alpha\in\Ult(V_\kappa,U)$.
\end{lemma}
\begin{proof}
  We prove this by induction on $\alpha$. Of course, this is trivial for $\alpha\leq\kappa$. Suppose that for every $x\in V_\alpha$ there is some $U$ for which $x\in\Ult(V_\kappa,U)$, and let $\gamma_x$ be the least such $\gamma$. We let $\delta$ be a regular cardinal such that for all $x\in V_\alpha$, $|\meas_\kappa(\gamma_x)|<\delta$.

  Let $j\colon V_\beta\to M$ be a $\delta$-supercompact embedding. Let $D=D(j,\sup j``\delta)$. As $\LP(\kappa)$ holds, $j_D\colon V_\beta\to M_D=\Ult(V_\beta,D)$ is an elementary embedding, and moreover, there is an elementary factor embedding $k\colon M_D\to M$ such that $j=k\circ j_D$ and $k([\id]_D) = \sup j``\delta$. Let us check that $\crit(k)\geq\delta$ and that $\meas_\kappa(\gamma_x)\in M_D$ for all $x\in V_\alpha$.

  First, we show that $j``\delta\in\rng(k)$. As $k$ is an elementary embedding, $\rng(k)$ is closed under definability in $M$. By \autoref{lemma:stationary-coding}, $j``\delta$ is definable from $\sup j``\delta$ in $M$. As $\sup j``\delta = k([\id]_D)$ is in the range of $k$, it follows that $j``\delta\in\rng(k)$ as wanted.

  We can now show that $k\restriction\delta=\id$, that is to say, $\crit(k)\geq\delta$. Note that $j``\delta=k``j_D``\delta$. As $\alpha\in j``\delta$ if and only if $\alpha=j(\beta)$ for some $\beta<\delta$, it follows that $\alpha=k(j_D(\beta))$ for some $\beta<\delta$. Therefore, $\alpha=k(\xi)$ for some $\xi$, and moreover, $\xi=j_D(\beta)$ for some $\beta<\delta$. It is now immediate that $j_D``\delta\in M$ and $k(j_D``\delta)=j``\delta$.

  As $j_D\restriction\delta$ is the unique isomorphism of $\delta$ and $j_D``\delta$, we have that $j_D\restriction\delta\in M_D$. Consequently, $k(j_D\restriction\delta)\restriction\delta=j\restriction\delta=(k\circ j_D)\restriction\delta$. On the other hand, the equality $(k\circ j_D)\restriction\delta = (k(j_D\restriction\delta)\circ k)\restriction\delta$ follows from the general fact on composition and application of elementary embeddings. Therefore, \[(k(j_D\restriction\delta)\circ k)\restriction\delta = k(j_D\restriction\delta)\restriction\delta.\] By cancelling out the $k(j_D\restriction\delta)$ parts we get that $k\restriction\delta=\id$ as wanted.

  Next, $\power(\delta)\in M_D$, as for every $S\subseteq\delta$, $j(S)\cap j_D``\delta=j``S$, and since $j_D\restriction\delta$ is the order isomorphism of $\delta$ and $j_D``\delta$, it is in $M_D$, so $S$ can be recovered from $j``S$. Therefore, $\delta^\delta\in M_D$ as well.

  Let $\vec U=\tup{U_\xi\mid\xi<\eta}$ be an enumeration of $\meas_\kappa(\gamma)$ for some $\gamma<\delta$. By the $\delta$-supercompactness of $j$, $\vec U\in M$, and since $k``M_D$ is an elementary submodel of $M$ and $\kappa,\gamma<\delta<\crit(k)$. Therefore, there is some $\vec W\in M_D$ such that $k(\vec W)=\vec U$. But, as $\crit(k)>\eta$, it must be that $k(\vec W)=k``\vec W$, and indeed that $k(W_\xi)=U_\xi$. So, $\meas_\kappa(\gamma)\in M_D$ as wanted.

  Note that every $f\in V_\kappa^\gamma$ can be decomposed into $f_1\colon\gamma\to\kappa$ and $f_2\colon\kappa\to V_\kappa$. Since $\delta^\delta\in M_D$, $(\kappa^\gamma)^{M_D}=\kappa^\gamma$, and since $V_{\kappa+1}\in M_D$, it means that $M_D$ computes the ultrapowers of $V_\kappa$ correctly. As $\meas_\kappa(\gamma)\in M_D$, it means that $x\in M_D$ for all $x\in V_\alpha$ for which $\delta>\gamma_x$. In particular, $V_\alpha\in M_D$ as well.

  Finally, as $\Ult(V_\kappa,D)$ is $V_{j_D(\kappa)}^{M_D}$, we have that $V_\alpha\in\Ult(V_\kappa,D)$. To complete the argument, we must show that every subset of $V_\alpha$ lies in an ultrapower as well. We use Solovay's capturing argument, applied to $M$ and $M_D$.

  Suppose that not every subset of $V_\alpha$ was in some ultrapower of $V_\kappa$. Then, $M$ would recognise this fact, being closed under $V_\delta$-sequences. Therefore, by elementarity, $M_D$ would have a counterexample as well, $Z\subseteq V_\alpha$. However, $k(Z)=Z$, so $M$ would recognise that $Z$ is a counterexample, despite computing $\Ult(V_\kappa,D)$ correctly.
\end{proof}

\begin{corollary}
  $V=L(\power^2(\Ord),V_{\kappa+1})=\HCD_{\kappa,V_{\kappa+1}}$.
\end{corollary}
\begin{proof}
  Note that $\kappa^\delta$ and $\meas_\kappa(\delta)$ are contained in $L(\power^2(\Ord))$ and $V_\kappa^\kappa\subseteq V_{\kappa+1}$. Therefore, as in the proof of \autoref{lemma:it-capturing},  $L(\power^2(\Ord),V_{\kappa+1})$ computes all the ultrapowers of $V_\kappa$ by $D\in\meas_\kappa(\delta)$ correctly, and therefore $V=L(\power^2(\Ord),V_{\kappa+1})$.

  To see that this is also equal to $\HCD_{\kappa,V_{\kappa+1}}$, first note that $V_{\kappa+1}$ is itself ordinal definable, so $\HCD_{\kappa,V_{\kappa+1}}=\HOD_{\meas_\kappa\cup V_{\kappa+1}\cup\{V_\kappa+1\}}$. Therefore, $V_{\kappa+1}\in\HCD_{\kappa,V_{\kappa+1}}$ and by \autoref{lemma:ord-hcd-capture} so is $\kappa^\delta$ for all $\delta\in\Ord$. So, by \autoref{lemma:it-capturing}, for every $x$ there is some $\delta$ and $U\in\meas_\kappa(\delta)$, and some $f\colon\delta\to\kappa$ and $g\colon\kappa\to V_\kappa$ such that $x=[g\circ f]_U$.
\end{proof}
\begin{proof}[Proof of \autoref{thm:spc+lp}]
  Under our assumptions, $V=\HCD_{\kappa,V_{\kappa+1}}$, which extends some $M=\HCD_{\kappa,\{a\}}\models\ZFC$. Let $c$ be a $V$-generic for $\Col(\omega,V_{\kappa+1})$, then in $V[c]$ every set in $V$ can be well-ordered, as the forcing is weakly homogeneous and ordinal definable.

  As $V$ surjects onto $V[c]$, this means that $\ZFC$ holds in $V[c]$. So by \autoref{fact:SVC} and $|V_{\kappa+1}^{<\omega}|=|V_{\kappa+1}|$, we have that $V_{\kappa+1}$ is a seed for $\SVC$ in $V$.
\end{proof}
\section{\texorpdfstring{$\SVC$}{SVC} in \texorpdfstring{$V$}{V} and large cardinals}
The dichotomy theorem tells us that in the presence of large cardinals, $\SVC$ fails significantly or else it holds in inner models. In this section, we will show that if $\SVC$ holds in $V$ itself, an extendible cardinal must already be above the rank of some seed, but this is not necessarily the case for a supercompact cardinal. We will do this by constructing a symmetric extension which preserves the supercompactness of a fixed cardinal. Towards that we prove a new lifting criterion for supercompact embeddings.
\subsection{Extendible cardinals and \texorpdfstring{$\SVC$}{SVC}}
\begin{definition}
  We say that $\delta$ is an \emph{extendible cardinal} if for every $\alpha\geq\delta$ there is some $\beta$ and a nontrivial elementary embedding $j\colon V_\alpha\to V_\beta$ with $\crit(j)=\delta$ and $j(\delta)>\alpha$.
\end{definition}
\begin{theorem}\label{thm:extnd-svc}
  Let $\delta$ be an extendible cardinal. Suppose that $\SVC$ holds, then for some $S\in V_\delta$, $\SVC(S)$ holds.
\end{theorem}
\begin{proof}
  Suppose that $\SVC(X)$ holds and let $\alpha$ be a $\beth$-fixed point such that $X\in V_\alpha$. Let $j\colon V_\alpha\to V_\beta$ be an elementary embedding with $\crit(j)=\delta$ and $j(\delta)>\alpha$. Since being a $\beth$-fixed point is a first-order property of $V_\alpha$, it must be that $\beta$ is a $\beth$-fixed point.
  By \autoref{lemma:strong-limit-svc}, $V_\beta\models\SVC(X)$. Therefore, $V_\beta\models\exists S(S\in V_{j(\delta)}\land\SVC(S))$. By elementarity, $V_\alpha\models\exists S(S\in V_\delta\land\SVC(S))$. As we may take $S=V_\gamma$ for some $\gamma<\delta$, this defines a function, mapping $\alpha$ to the least $\gamma$ for which $V_\alpha\models\SVC(V_\gamma)$. By the pigeonhole principle, on a proper class of ordinals $\gamma$ is fixed, and so $\SVC(V_\gamma)$ holds for some $\gamma<\delta$.
\end{proof}

While this paper is concerned with models of $\ZF$ where the Axiom of Choice fails, throughout the section our ground model will be a model of $\ZFC$. In particular, we will not worry about the type of supercompactness in $V$.
\subsection{Symmetric extensions and lifting arguments}\label{subsect:sym-ext}
Symmetric extensions allow us to extend a model of $\ZFC$ and violate the Axiom of Choice. This is done by identifying a particular model between the ground model and the generic extension. For a complete study, see \cite{KS:Towards}. We include a brief summary of the technique below.

Let $\PP$ be a notion of forcing. If $\sG$ is a group of automorphisms of $\PP$,\footnote{Often we consider a group acting on $\PP$, but we will ignore this distinction in favour of readability.} its action extends to $\PP$-names by recursion, \[\pi\dot x = \{\tup{\pi p,\pi\dot y}\mid \tup{p,\dot y}\in\dot x\}.\]

\begin{lemma}[The Symmetry Lemma]
  Let $p\in\PP$, let $\dot x$ be a $\PP$-name, let $\varphi$ be a formula in the language of forcing, and let $\pi$ be an automorphism of $\PP$. Then, \[\pushQED{\qed}p\forces\varphi(\dot x)\iff\pi p\forces\varphi(\pi\dot x).\qedhere\popQED\]
\end{lemma}

Let $\sF$ be a \emph{filter of groups on $\sG$}, that is a family of subgroups of $\sG$ closed under supergroups, finite intersections, and not containing the trivial group. We say that a filter of groups is \emph{normal} if it is closed under conjugations by the elements of $\sG$.\footnote{This is an unfortunate clash between normal filters in this definition and in the sense of large cardinals. In this work, we will only discuss ultrafilters and refer to them as \emph{measures} to avoid any possible confusion.}

A \emph{symmetric system} is a triple $\cS=\tup{\PP,\sG,\sF}$ where $\PP$ is a notion of forcing, $\sG$ is a group of automorphisms of $\PP$, and $\sF$ is a normal filter of subgroups of $\sG$. Let us fix such $\cS$.

We call a $\PP$-name, $\dot x$, a \emph{$\cS$-symmetric} name if $\sym_\sG(\dot x)=\{\pi\in\sG\mid\pi\dot x=\dot x\}\in\sF$. We let $\HS_\cS$ be the class of all \emph{hereditarily $\cS$-symmetric}, also called $\cS$-names. This is defined recursively and we denote by $\HS_{\cS,\alpha}$ the class of $\cS$-names constructed in the first $\alpha$ stages. When context is clear we will omit all subscripts.

The class $\HS_\cS$ has a forcing relation associated with it, given by the relativisation of quantifiers to the class $\cS$. Namely, $p\forces_\cS\varphi(\dot x)$ means that $\dot x\in\cS$ and that every quantifier is limited to the names in $\cS$ as well. This forcing relation is often not full (that is, $p\forces_\cS\exists x\varphi(x)$ need not imply that there is some $\dot x\in\HS$ such that $p\forces_\cS\varphi(\dot x)$). When $\pi\in\sG$ we have that $p\forces_\cS\varphi(\dot x)$ if and only if $\pi p\forces_\cS\varphi(\pi\dot x)$.

\begin{theorem}
  Let $\cS$ be a symmetric system, and let $G\subseteq\PP$ be a $V$-generic filter. Then $\HS^G=\{\dot x^G\mid\dot x\in\HS\}$ satisfies $\ZF$ and contains $V$.
\end{theorem}
The model $\HS^G$ is called a \emph{symmetric extension} of $V$. If $V$ satisfies $\ZFC$, then its symmetric extensions satisfy $\SVC$. Conversely, any model of $\SVC$ is a symmetric extension of some model of $\ZFC$, as was shown by Usuba \cite{Usuba:2021}.

Much like with forcing, we can take the products of symmetric systems, where the automorphism groups are acting pointwise, and the filters are generated by the products of groups in each filter.\footnote{This can be extended to iterations, but we will not need that for our purposes.}

When extending the ground model we want to make sure that large cardinals remain large in the extension, which in this case requires us to find suitable elementary embeddings. In the case of generic extensions we have Silver's criterion to determine exactly when an embedding can be lifted, and moreover when is this lifting can be seen by the extension.

Hayut and Karagila devised such a condition in \cite{HK:Critical}. This was further simplified by Karagila and Yuan in \cite[Theorem~2.3]{KY:2026}. We rephrase this simplified lifting theorem to make its application neater in our case.

\begin{theorem}[Simplified Lifting Theorem]\label{thm:lifting}
  Suppose that $j\colon V\to M$ is an elementary embedding and $\cS=\tup{\PP,\sG,\sF}$ is a symmetric system such that $j(\cS)\cong\cS\times\cT$, where $\cT=\tup{\RR,\sH,\sK}$ is a symmetric system in $M$ with the following properties:
  \begin{enumerate}
  \item There exists an $M$-generic filter for $\RR$ in $V$.
  \item For every $p\in\PP$, $j(p)=\tup{p,\1_\RR}$.
  \item For every $\sigma\in\sG$, $j(\sigma)=\tup{\sigma,\id}$.
  \end{enumerate}
  Then we can lift $j$ amenably to the symmetric extensions of $V$ and $M$ by $\cS$ and $j(\cS)$ respectively.
\end{theorem}

In \cite{KY:2026} the theorem required that for all $K\in j(\sF)$, $\{\sigma\mid j(\sigma)\in K\}\in\sF$. This, in the presence of the other conditions, is equivalent to $j``(K\restriction\sG)\subseteq K$,\footnote{Here $K\restriction\sG$ is the projection of $K$ onto a subgroup of $\sG$.} which then nearly trivialises when we require that $j(\cS)\cong\cS\times\cT$: If $K\in j(\sF)$, then we may assume that $K=K_0\times K_1$ where $K_0\in\sF$ and $K_1\in\sK$. So, the requirement is that $j``K_0\subseteq K$. However, $j``K_0=K_0\times\{\id\}$, which easily is a subgroup of $K$.

Lifting the embeddings from the ground model is not enough, though. We need to guarantee that supercompact embeddings remain supercompact, namely that they have sufficiently closed target models.

For this we prove a new criterion to preserve the closure of the target models, under the conditions of the Simplified Lifting Theorem. First, we use the concept of a partial interpretation to transform a name in $\HS_{j(\cS)}^M$ to a name in $\HS_\cS$ using the $M$-generic filter, $H$. This was defined in \cite{HK:Critical} under more general conditions, but we can simplify it in our case. We define for $\dot x\in\HS_{j(\cS)}^M$ by recursion,
\[\dot x^H = \{\tup{p,\dot y^H}\mid\tup{\tup{p,r},\dot y}\in\dot x,r\in H\}.\]
So, for example, under the conditions of the lifting theorem above, $(\dot x^H)^G=\dot x^{G\times H}$, where $G$ is a $V$-generic filter for $\PP$.
\begin{proposition}\label{prop:sym-p-int}
If $\tup{\pi,\sigma}\in j(\sG)$, then $(\tup{\pi,\sigma}\dot x)^H=\pi(\dot x^{\sigma``H})$.
\end{proposition}
\begin{proof}
  By induction on the rank of $\dot x$,
  \begin{align*}
    (\tup{\pi,\sigma}\dot x)^H &= (\{\tup{\tup{\pi p,\sigma r},\tup{\pi,\sigma}\dot y}\mid\tup{\tup{p,r},\dot y}\in\dot x\})^H\\
                               &= \{\tup{\pi p,\pi(\dot y^{\sigma``H})}\mid\tup{\tup{p,r},\dot y}\in\dot x,\sigma r\in H\}\\
                               &= \{\tup{\pi p,\pi(\dot y^{\sigma``H})}\mid\tup{\tup{p,r},\dot y}\in\dot x,r\in\sigma``H\}\\
                               &= \pi\{\tup{p,\dot y^{\sigma``H}}\mid\tup{\tup{p,r},\dot y}\in\dot x,r\in\sigma``H\}\\
                               &= \pi\dot x^{\sigma``H}.\qedhere
  \end{align*}
\end{proof}

\begin{theorem}\label{thm:spc-lifting}
  Suppose that the conditions of the Simplified Lifting Theorem hold. Moreover, assume that $j$ is a $\lambda$-supercompact embedding where $\lambda=|V_\lambda|$. We can lift $j$ to a $\lambda$-supercompact embedding if the following conditions hold:
  \begin{enumerate}
  \item $|\PP|<\lambda^+$,
  \item $\sK$ is a $\lambda^+$-complete filter of subgroups.
  \end{enumerate}
\end{theorem}
\begin{proof}
  Following the assumption of the theorem, let $W$ and $N$ denote the symmetric extensions of $V$ and $M$. Suppose $f\colon W_\lambda\to N$ is a function in $W$, we want to show that $f\in N$. Let $\dot f\in\HS_\cS$ be a name for $f$. By the fact that $|\PP|<\lambda^+$, we may assume that $|\dot f|=\lambda$. Moreover, by refining $\dot f$ if necessary, we may assume that if $\tup{\dot x,\dot y'}^\bullet$ appears in $\dot f$, then $\dot y'=\dot y^H$ for some $\dot y\in\HS_{j(\cS)}^M$.

  First, note that $H_{\lambda^+}\subseteq M$, so if $\dot x$ is a $\PP$-name such that $|\tcl(\dot x)|\leq\lambda$, then $\dot x\in M$. In particular, $W_\lambda$, whose canonical name is $\HS_{\cS,\lambda}^\bullet$ is in $M$. Moreover, every $\PP$-name which lies in $M$ has a natural extension to a $j(\PP)$-name. So, if $\dot x\in\HS_{\cS,\lambda}$, then $\dot x\in\HS_{j(\cS)}^M$ and $\dot x^H=\dot x$. We will therefore ignore the distinction for $\dot x\in\HS_{\cS,\lambda}$ between its $\PP$-name and $j(\PP)$-name.

  It follows now that the $j(\PP)$-name $\dot f_*=\{\tup{\tup{p,\1},\tup{\dot x,\dot y}^\bullet}\mid\tup{p,\tup{\dot x,\dot y^H}^\bullet}\in\dot f\}$ is in $M$. Let us verify that this name is a $j(\cS)$-name. For each $\dot y$ which appears in the range of $\dot f_*$ let $K_{1,\dot y}$ denote the projection of $\sym_{j(\sG)}(\dot y)$ onto $\cT$. Let $K_1=\bigcap K_{1,\dot y}$, ranging over all the $\dot y$.\footnote{Note that for all $\dot x$ which appear in the domain of $\dot g$, said projection is $\sH$ itself and can be ignored.} As $\sK$ is $\lambda^+$-complete, we have that $K_1\in\sK$ and so $\sym(\dot f)\times K_1\in j(\sF)$.

  We claim that if $\tup{\pi,\sigma}\in\sym(\dot f)\times K_1$, then $\tup{\pi,\sigma}\dot f_*=\dot f_*$. As all the names which appear in $\dot f_*$ are in $\HS_{j(\cS)}^M$, this would complete the proof. As $\tup{\pi,\sigma}=\tup{\id,\sigma}\circ\tup{\pi,\id}$, let us tackle the action of $\pi$ and $\sigma$ separately.

  By \autoref{prop:sym-p-int}, $\tup{\pi,\id}\dot f_*^H=\pi(\dot f_*^H)$, but as it is immediate from the definition that $\dot f_*^H=\dot f$ and $\pi\in\sym(\dot f)$, we get that $\tup{\pi,\id}\dot f_*=\dot f_*$.

  Let us inspect the application of $\tup{\id,\sigma}$ to a single ordered pair, $\tup{\tup{p,\1_\RR},\tup{\dot x,\dot y}^\bullet}$, in $\dot f_*$. For the condition, $\tup{\id,\sigma}(\tup{p,\1_\RR})=\tup{p,\1_\RR}$, and so it is not changed at all. For $\tup{\dot x,\dot y}^\bullet$, as $\dot x\in\HS_{\cS,\lambda}$, we immediately have that $\tup{\id,\sigma}\dot x=\dot x$, and since $\sigma\in K_1$, $\tup{\id,\sigma}\dot y=\dot y$ as well. Therefore, $\tup{\id,\sigma}\dot f_*=\dot f_*$ as wanted.
\end{proof}
\subsection{Supercompactness does not reflect \texorpdfstring{$\SVC$}{SVC}}
\begin{theorem}\label{thm:spc+svc}
  Let $\kappa$ be the least supercompact cardinal. There is a symmetric extension, $W$, in which $\kappa$ is supercompact, $\SVC(W_\kappa)$ holds, but $\SVC(W_\alpha)$ fails for all $\alpha<\kappa$.
\end{theorem}
\begin{proof}
  We let $A$ be the set $\{\alpha<\kappa\mid\forall\delta<\alpha,\delta\text{ is not }\alpha\text{-supercompact}\}$. For each $\alpha\in A$ we let $\cS_\alpha=\tup{\PP_\alpha,\sG_\alpha,\sF_\alpha}$, where
  \begin{enumerate}
  \item $\PP_\alpha=\Add(\alpha^+,\alpha^+)$. Namely, $p\colon\alpha^+\times\alpha^+\to 2$ with $|p|<\alpha^+$.
  \item $\sG_\alpha=S_{<\alpha^+}$ is the permutations of $\alpha^+$ which only move fewer than $\alpha^+$ points, with the action on $\PP_\alpha$ given by $\pi p(\pi\eta,\xi)=p(\eta,\xi)$.
  \item $\sF_\alpha$ is generated by $\fix(E)=\{\pi\in\sG_\alpha\mid\pi\restriction E=\id\}$ for $E\in\power_{\alpha^+}(\alpha^+)$.
  \end{enumerate}

  We let $\cS$ be the Easton support product of $\cS_\alpha$. Namely, $\PP$ is the Easton support product of $\PP_\alpha$; $\sG$ is the Easton support product of the groups, acting pointwise on $\PP$; and $\sF$ is generated by sequences of groups from the filter that on an Easton-small set are not $\sG_\alpha$ itself. Fixing a $V$-generic filter $G\subseteq\PP$, let $W$ denote the symmetric extension $\HS_\cS^G$.

  Let $\lambda>\kappa$ be such that $|V_\lambda|=\lambda$, $\cf(\lambda)>\kappa$, and $2^\lambda=\lambda^+$.\footnote{The final condition holds for any singular cardinal above $\kappa$ as a consequence of $\SCH$ holding above supercompact cardinals.} Let $j\colon V\to M$ be an ultrapower embedding with a supercompact measure on $\power_\kappa(\lambda)$, so that $M^\lambda\subseteq M$. Let us verify that the conditions for \autoref{thm:lifting} and \autoref{thm:spc-lifting} hold, so that $j$ lifts to a $\lambda$-supercompact embedding in $W$.

  First, note that $j(\PP)=\prod_{\alpha\in j(A)}\cS_\alpha$, where $\cS_\alpha$ are similarly defined as above, which allows us to define, in $M$, the Easton support product \[\cT=\tup{\RR,\sH,\sK}=\prod_{\alpha\in j(A)\setminus A}\cS_\alpha.\]

  Easily, $\PP\subseteq V_\kappa$, so $|\PP|=\kappa<\lambda^+$. Since $M^\lambda\subseteq M$ and $\RR$ is $\lambda^+$-closed, and since $j$ is an ultrapower embedding, $\RR$ has an $M$-generic filter in $V$.

  Next, as Easton support is bounded below $\kappa$, conditions (2) and (3) in \autoref{thm:lifting} hold automatically. Finally, as each $\sF_\alpha$ is $\alpha^+$-complete, it must be that $\sK$ is at least $\lambda^+$-complete. And so, $\kappa$ remains $\lambda$-supercompact in $W$ for all $\lambda$.

  For each $\alpha\in A$ we let $\dot x^\alpha_\eta=\{\tup{p,\check\xi}\mid p(\eta,\xi)=1\}$ and let $\dot X^\alpha=\{\dot x_\eta^\alpha\mid\alpha<\alpha^+\}^\bullet$. It is not hard to see that $\pi\dot x^\alpha_\eta=\dot x^\alpha_{\pi\eta}$ and so each of those is a $\cS_\alpha$-name and therefore $\sym(\dot X^\alpha)=\sG_\alpha$. We omit the dots to indicate the interpretation of the name in our symmetric extension.

  To complete the proof we need to show that $\SVC(W_\alpha)$ fails for all $\alpha<\kappa$, while $\SVC(W_\kappa)$ holds. Let us fix $\alpha<\kappa$, and let $\dot f\in\HS_\cS$ be a name for a function from $W_\alpha\times\beta\to X^\delta$ where $\delta$ is some ordinal such that $|\HS_{\alpha,\cS}|<\delta$. Let $\vec e$ be a sequence of sets such that $\tup{\fix(e_\alpha)\mid\alpha<\kappa}\subseteq\sym(\dot f)$. Let $p\in\PP$ be a condition such that $p\forces_\cS\dot f\colon\HS_{\alpha,\cS}^\bullet\times\check\eta\to\dot X^\delta$, and without loss of generality $\supp p_\alpha\subseteq e_\alpha$ for all $\alpha<\kappa$.

  Let $\xi<\delta^+$ be such that $\xi\notin e_\delta$, and let $q\leq p$ be a condition such that for some $\dot x\in\HS_{\alpha,\cS}$ and $\gamma<\beta$, $q\forces_\cS\dot f(\dot x,\check\gamma)=\dot x^\delta_\xi$. Pick some $\xi'\notin e_\delta\cup\{\xi\}\cup\{\supp q_\delta\}$ and let $\pi$ be the automorphism such that $\pi_\delta$ is the $2$-cycle $(\xi\ \xi')$, and otherwise the identity. As $\delta>\alpha$, $\pi\dot x=\dot x$ and $\pi\dot f=\dot f$. Therefore, $\pi q\forces_\cS\dot f(\dot x,\check\gamma)=\dot x^\delta_{\xi'}\neq\dot x^\delta_\xi$. This is impossible, since it $q$ and $\pi q$ are compatible, so any common extension must force that $\dot f$ is not a function. Therefore, it must be that $p$ forced that $\rng f\subseteq\{x^\delta_\xi\mid\xi\in e_\delta\}$, and so it is not surjective. Therefore, $W\models\lnot\SVC(W_\alpha)$, for any $\alpha<\kappa$.

  Finally, for any sequence, $e$, such that $\tup{\fix(e_\alpha)\mid\alpha<\kappa}\in\sG$, we let $\dot f_e$ be the name $\{\tup{\check\alpha,\check\xi,\dot x^\alpha_\xi}^\bullet\mid\xi\in e_\alpha\}^\bullet$. Using the terminology and notation of \cite[\S7]{KS:Towards}, the collection of all such $\dot f_e$ forms a \emph{respect basis} for $\cS$. Moreover, as $e_\alpha=\varnothing$ for a tail of $\alpha<\kappa$, $\dot f_e\in V_\kappa$. It follows from this that the quotient $\PP/\cS$ is contained in $W_\kappa$, and by \cite[Theorem~8.1]{KS:Towards}, $W=V(\PP/\cS)$. Therefore, $\SVC(W_\kappa)$ holds.
\end{proof}
\begin{remark}
  For any fixed $\alpha<\kappa$ it is possible for $\DC_\alpha$ to hold in the model constructed above as well. To see this, simply restrict $A$ such that $\min A>\alpha$, we have that $\PP$ is $\alpha^+$-closed and $\sF$ is $\alpha^+$-complete. So by \cite{Karagila:2019}, $\DC_\alpha$ holds in the model.
\end{remark}
\begin{remark}
  The above theorem can be extended to any supercompact cardinal. Fixing a function $f\colon\kappa\to\kappa$ such that for any $\lambda$ there is a supercompact measure on $\power_\kappa(\lambda)$ such that $j_U(f)(\kappa)>\lambda$, we let $A$ be the set of closure points of $f$ and proceed as before.
\end{remark}
\begin{theorem}
  In the model constructed in \autoref{thm:spc+svc}, $\LP(\kappa)$ fails.
\end{theorem}
\begin{proof}
  We use the same notation from the previous proof. If $p\in\PP$, we write $p_\alpha$ for its projection onto $\PP_\alpha$, and similarly for $\pi\in\sG$.

  First, note that the sequence $\tup{X^\alpha\mid\alpha<\kappa}$ is in $W$, as any automorphism in $\sG$ must preserve the canonical name for this sequence. Suppose now that $\dot f\in\HS_\cS$ such that $p\forces_\cS``\dot f(\check\alpha)\in\dot X^\alpha$ for unboundedly many $\alpha$''.

  Let $\alpha$ be large enough such that $\supp(p)\subseteq\alpha$ and that $\supp(\sym(\dot f))\subseteq\alpha$.\footnote{Here, $\supp(\sym(\dot f))$ would be the coordinates on which the projection is not $\sG_\alpha$.} Suppose that $q\leq p$ and $q\forces_\cS\dot f(\check\alpha)=\dot x^\alpha_\eta$ for some $\eta<\alpha^+$.

  Let $\eta'<\alpha^+$ such that $\eta'\notin\supp(q_\alpha)$. Let $\pi\in\sG$ be such that $\pi_\beta=\id$ for all $\beta\neq\alpha$ and $\pi_\alpha$ is the $2$-cycle $(\eta\ \eta')$. Then $\pi\in\sym(\dot f)$, $\pi p=p$, and $\pi q$ is compatible with $q$. However, $\pi q\forces_\cS\dot x^\alpha_{\eta'}=\pi\dot x^\alpha_\eta=\pi\dot f(\pi\check\alpha)=\dot f(\check\alpha)\neq\dot x^\alpha_\eta$. This is impossible, as any extension of $q$ and $\pi q$ would then force that $\dot f$ is not a function.

  Therefore, in $W$, there is no unbounded $I\subseteq\kappa$ for which $\prod_{\alpha\in I}X^\alpha\neq\varnothing$. In particular, for any $U\in\meas_\kappa(\kappa)$, $U$-$\AC_\kappa$ must fail. As $\kappa$ is a supercompact cardinal, it is measurable so $\meas_\kappa(\kappa)$ is non-empty. Therefore, $\LP(\kappa)$ fails.
\end{proof}
\section{Conclusions}
\begin{question}
  Suppose that an extendible cardinal exists and $\ACH$ holds. Does $\SVC$ hold in $V$?
\end{question}
The dichotomy theorem shows that $\SVC$ holds in $L(\power(\Ord))$. But what about $V$ itself? Can we somehow use the extendible cardinal to pull $\SVC$ from $L(\power(\Ord))$ to $V$ itself? If so, the dichotomy turns into the following statement: In $\delta$ is an extendible cardinal, then either all regular cardinals above $\delta$ are $\omega$-J\'onsson, or $\SVC(V_\alpha)$ holds for some $\alpha<\delta$.

Woodin proved that if $\kappa$ is a supercompact cardinal, then the finite support $\prod_{\alpha<\kappa}\Col(\omega,V_\alpha)$ forces $\DC$ to hold in the universe \cite[Lemma~225]{Woodin:SEM1}. Hayut and Karagila used this as a criterion for the failure of $\AC$ with supercompactness: the failure is trivial, when Woodin's forcing actually restores $\AC$ and not just $\DC$. This implies that $\SVC(V_\kappa)$ holds, when the failure is trivial.\footnote{This was informed by the mistaken idea that $\SVC$ will reflect below the supercompact, in which case the failure would be trivial if and only if $\SVC(V_\alpha)$ holds for some $\alpha<\kappa$. As we see in this paper, this is not quite true.}

\begin{question}[Hayut--Karagila {\cite{HK:Critical}}]
  Is a nontrivial failure of choice with supercompactness consistent relative to $\ZFC$ large cardinals?
\end{question}

One trivial upper bound is a super-Reinhardt cardinal (and less trivially, one Reinhardt would also suffice). But one would hope that this can be done starting with a single supercompact in a model of $\ZFC$. While we do not have an answer, the work in this paper shows that this question cannot be solved using a single symmetric extension when starting from a model of $\ZFC$. Indeed, if this is possible at all, this will require a class-sized symmetric extension (which or may not be an $\Ord$-length iteration of set-sized symmetric extensions).

Indeed, the original of the Bristol model \cite{Bristol1,Bristol2}, was trying to understand the Axiom of Choice Conjecture, and in this model, $L(\power^\alpha(\Ord))\models\SVC$ for all $\alpha$, even though $\SVC$ fails. It may be that one can find a similar construction which preserves the supercompactness, however the Bristol model construction cannot be carried above with a supercompact present in the universe.

We finish the paper by discussion the obvious connections between the Axiom of Choice Dichotomy and the HOD Dichotomy. Recall that $\delta$ is \emph{$\omega$-strongly measurable in $\HOD$} if for some $\gamma<\delta$ we have: (1) $\HOD\models|\gamma|=\gamma$ and $2^\gamma<\kappa$; (2) there is no partition of $S^\kappa_\omega$ into $\gamma$ stationary sets in $\HOD$.

The HOD Dichotomy says that in the presence of a strongly compact cardinal, either all sufficiently large singular cardinals are singular in $\HOD$ and their successor is computed correctly, or else all sufficiently large regular cardinals are $\omega$-strongly measurable in $\HOD$.

We can strengthen \autoref{thm:spc-dich} by replacing (1) with the statement that there exists a cardinal $\lambda$ such that for all regular cardinals $\delta>\lambda$, any $f\colon[\delta]^\omega\to\lambda$ omits a colour on an $\omega$-club subset of $\delta$. Easily, any such $\delta$ is $\omega$-strongly measurable in $\HOD$, and so the HOD Hypothesis must fail in this case.

\providecommand{\bysame}{\leavevmode\hbox to3em{\hrulefill}\thinspace}
\providecommand{\MR}{\relax\ifhmode\unskip\space\fi MR }
\providecommand{\MRhref}[2]{%
  \href{http://www.ams.org/mathscinet-getitem?mr=#1}{#2}
}
\providecommand{\href}[2]{#2}

\end{document}